\documentstyle[12pt,fleqn]{article}

\begin{document}

\title{Planar Complex Numbers in Even $n$ Dimensions}

\author{Silviu Olariu
\thanks{e-mail: olariu@ifin.nipne.ro}\\
Institute of Physics and Nuclear Engineering,\\
Department of Fundamental Experimental Physics\\
76900 Magurele, P.O. Box MG-6, Bucharest, Romania}

\date{4 August 2000}

\maketitle

\abstract

Planar commutative n-complex numbers of the form
$u=x_0+h_1x_1+h_2x_2+\cdots+h_{n-1}x_{n-1}$ are introduced in an even number n
of dimensions, the variables $x_0,...,x_{n-1}$ being real numbers.  The planar
n-complex numbers can be described by the modulus $d$, by the amplitude $\rho$,
by $n/2$ azimuthal angles $\phi_k$, and by n/2-1 planar angles $\psi_{k-1}$.
The exponential function of a planar n-complex number can be expanded in terms
of the planar n-dimensional cosexponential functions $f_{nk}, k=0,1,...,n-1$,
and expressions are given for $f_{nk}$.  Exponential and trigonometric forms
are obtained for the planar n-complex numbers.  The planar n-complex functions
defined by series of powers are analytic, and the partial derivatives of the
components of the planar n-complex functions are closely related.  The
integrals of planar n-complex functions are independent of path in regions
where the functions are regular.  The fact that the exponential form of a
planar n-complex numbers depends on the cyclic variables $\phi_k$ leads to the
concept of pole and residue for integrals on closed paths.  The polynomials of
planar n-complex variables can always be written as products of linear factors,
although the factorization may not be unique.

\endabstract

\section{Introduction}

A regular, two-dimensional complex number $x+iy$ 
can be represented geometrically by the modulus $\rho=(x^2+y^2)^{1/2}$ and 
by the polar angle $\theta=\arctan(y/x)$. The modulus $\rho$ is multiplicative
and the polar angle $\theta$ is additive upon the multiplication of ordinary 
complex numbers.

The quaternions of Hamilton are a system of hypercomplex numbers
defined in four dimensions, the
multiplication being a noncommutative operation, \cite{1} 
and many other hypercomplex systems are
possible, \cite{1}-\cite{2b} but these hypercomplex systems 
do not have all the required properties of regular, 
two-dimensional complex numbers which rendered possible the development of the 
theory of functions of a complex variable.

A system of hypercomplex numbers in an even number of dimensions $n$
is described in this work, 
for which the multiplication is both associative and commutative, and which is 
rich enough in properties so that an exponential form exists and the concepts
of analytic n-complex 
function,  contour integration and residue can be defined.
The n-complex numbers introduced in this work have 
the form $u=x_0+h_1x_1+h_2x_2+\cdots+h_{n-1}x_{n-1}$, the variables 
$x_0,...,x_{n-1}$ being real numbers. The multiplication rules for the complex
units $h_1,...,h_{n-1}$ are 
$h_j h_k =h_{j+k}$ if $0\leq j+k\leq n-1$, and $h_jh_k=-h_{j+k-n}$ if
$n\leq j+k\leq 2n-2$, where $h_0=1$.
The product of two n-complex numbers is equal to zero if both numbers are
equal to zero, or if the numbers belong to certain n-dimensional hyperplanes
described further in this work. 

If the n-complex number $u=x_0+h_1x_1+h_2x_2+\cdots+h_{n-1}x_{n-1}$ is
represented by the point $A$ of coordinates $x_0,x_1,...,x_{n-1}$, 
the position of the point $A$ can be described, in an even number 
of dimensions, by the modulus $d=(x_0^2+x_1^2+\cdots+x_{n-1}^2)^{1/2}$, 
by $n/2$ azimuthal angles $\phi_k$ and by $n/2-1$
planar angles $\psi_{k-1}$. 
An amplitude $\rho$ can be defined 
as $\rho^n=\rho_1^2\cdots\rho_{n/2}^2$,
where $\rho_k$ are radii in orthogonal two-dimensional planes defined further
in this work. The amplitude $\rho$, the radii
$\rho_k$ and the variables $\tan\psi_{k-1}$ are multiplicative, and the
azimuthal angles $\phi_k$ are 
additive upon the multiplication of n-complex numbers.
Because the description of the position of the point $A$ requires, in addition
to the azimuthal angles $\phi_k$, only the planar angles $\psi_{k-1}$, 
the hypercomplex numbers studied in this work will be called planar
n-complex number, to distinguish them from the polar n-complex numbers
introduced in \cite{1a}, which in an even number of dimensions required two
polar angles, and in an odd number of dimensions required one polar angle.

The exponential function of an n-complex number can be expanded in terms of
the planar n-dimensional cosexponential functions
$f_{nk}(y)=\sum_{p=0}^\infty (-1)^p y^{k+pn}/(k+pn)!, k=0,1,...,n-1$.
It is shown that
$f_{nk}(y)=(1/n)\sum_{l=1}^{n}$
$\exp\left\{y\cos\left(\pi (2l-1)/n\right)\right\}$\\
$\cos\left\{y\sin\left(\pi (2l-1)/n\right)-\pi k(2l-1)/n\right\}, $
$k=0,1,...,n-1.$
Addition theorems and other relations
are obtained for the planar n-dimensional cosexponential functions.

The exponential form of an n-complex number, which can be defined for
all $x_0,...,x_{n-1}$, is
$u=\rho \exp\left\{\sum_{p=1}^{n-1}h_p\left[
-(2/n)\sum_{k=2}^{n/2}
\cos\left(\pi (2k-1)p/n\right)\ln\tan\psi_{k-1}
\right]\right\}$
$\exp\left(\sum_{k=1}^{n/2}\tilde e_k\phi_k\right)$,
where $\tilde e_k=(2/n)\sum_{p=1}^{n-1}h_p\sin(\pi (2k-1)p/n)$. 
A trigonometric form also exists for an n-complex number,
$u=d\left(n/2\right)^{1/2}$
$\left(1+1/\tan^2\psi_1+1/\tan^2\psi_2+\cdots
+1/\tan^2\psi_{n/2-1}\right)^{-1/2}$\\
$\left(e_1+\sum_{k=2}^{n/2}e_k/\tan\psi_{k-1}\right)$
$\exp\left(\sum_{k=1}^{n/2}\tilde e_k\phi_k\right)$.

Expressions are given for the elementary functions of n-complex variable.
The functions $f(u)$ of n-complex variable which are defined by power series
have derivatives independent of the direction of approach to the point under
consideration. If the n-complex function $f(u)$ 
of the n-complex variable $u$ is written in terms of 
the real functions $P_k(x_0,...,x_{n-1}), k=0,...,n-1$, then
relations of equality  
exist between partial derivatives of the functions $P_k$. 
The integral $\int_A^B f(u) du$ of an n-complex
function between two points $A,B$ is independent of the path connecting $A,B$,
in regions where $f$ is regular.
If $f(u)$ is an analytic n-complex function, then 
$\oint_\Gamma f(u)du/(u-u_0)$
$=2\pi f(u_0)\sum_{k=1}^{n/2}\tilde e_k$ 
$\;{\rm int}(u_{0\xi_k\eta_k},\Gamma_{\xi_k\eta_k})$,
where the functional ${\rm int}$ takes the values 0 or 1 depending on the
relation between $u_{0\xi_k\eta_k}$ and $\Gamma_{\xi_k\eta_k}$, which are
respectively the projections of the point $u_0$ and of 
the loop $\Gamma$ on the plane defined by the orthogonal axes $\xi_k$ and
$\eta_k$, as expained further in this work.

A polynomial $u^m+a_1 u^{m-1}+\cdots+a_{m-1} u +a_m  $ can always
be written as a product of linear factors, although the factorization may not
be unique. 

This paper belongs to a series of studies on commutative complex numbers in $n$
dimensions. \cite{2c}
For $n=2$, the n-complex numbers discussed in this paper become the usual
2-dimensional complex numbers $x+iy$.
A detailed analysis for $n=4$ and $n=6$  of the planar
n-complex numbers can be found in the corresponding
studies mentioned in Ref. \cite{2c}.

\section{Operations with planar n-complex numbers}

A hypercomplex number in $n$ dimensions is determined by its $n$ components
$(x_0,x_1,...,x_{n-1})$. The planar n-complex numbers and
their operations discussed in this work can be represented 
by  writing the n-complex number $(x_0,x_1,...,x_{n-1})$ as  
$u=x_0+h_1x_1+h_2x_2+\cdots+h_{n-1}x_{n-1}$, where 
$h_1, h_2, \cdots, h_{n-1}$ are bases for which the multiplication rules are 
\begin{equation}
h_j h_k =(-1)^{[(j+k)/n]}h_l ,\:l=j+k-n[(i+k)/n],
\label{1}
\end{equation}
for $ j,k,l=0,1,..., n-1$,
where $h_0=1$.
In Eq. (\ref{1}), $[(j+k)/n]$ denotes the integer part of $(j+k)/n$, 
the integer part being defined as $[a]\leq a<[a]+1$, so that
$0\leq j+k-n[(j+k)/n]\leq n-1$. 
In this work, brackets larger than the regular brackets
$[\;]$ do not have the meaning of integer part.
The significance of the composition laws in Eq.
(\ref{1}) can be understood by representing the bases $h_j, h_k$ by points on a
circle at the angles $\alpha_j=\pi j/n,\alpha_k=\pi k/n$, as shown in Fig. 1,
and the product $h_j h_k$ by the point of the circle at the angle 
$\pi (j+k)/n$. If $\pi\leq \pi (j+k)/n<2\pi$, the point is opposite to the
basis $h_l$ of angle $\alpha_l=\pi (j+k)/n-\pi$.

In an odd number of dimensions $n$, a transformation of coordinates according
to
\begin{equation}
x_{2l}=x^\prime_l, x_{2m-1}=-x^\prime_{(n-1)/2+m}, 
\label{1a}
\end{equation}
and of the bases according to
\begin{equation}
h_{2l}=h^\prime_l, h_{2m-1}=-h^\prime_{(n-1)/2+m}, 
\label{1b}
\end{equation}
where $l=0,...,(n-1)/2, \; m=1,...,(n-1)/2$,
leaves the expression of an n-complex number unchanged,
\begin{equation}
\sum_{k=0}^{n-1}h_k x_k=\sum_{k=0}^{n-1}h^\prime_k x^\prime_k,
\label{1c}
\end{equation}
and the products of the bases $h^\prime_k$ are
\begin{equation}
h^\prime_j h^\prime_k =h^\prime_l ,\:l=j+k-n[(j+k)/n],
\label{1d}
\end{equation}
for $j,k,l=0,1,..., n-1$. 
Thus, the n-complex numbers with the rules (\ref{1}) are equivalent in an
odd number of dimensions to the n-complex numbers described in \cite{1a}.
Therefore, in this work it will be supposed that $n$ is an even number, unless
otherwise stated.

Two n-complex numbers 
$u=x_0+h_1x_1+h_2x_2+\cdots+h_{n-1}x_{n-1}$,
$u^\prime=x^\prime_0+h_1x^\prime_1+h_2x^\prime_2+\cdots+h_{n-1}x^\prime_{n-1}$ 
are equal if and only if $x_j=x^\prime_j, j=0,1,...,n-1$.
The sum of the n-complex numbers $u$
and
$u^\prime$ 
is
\begin{equation}
u+u^\prime=x_0+x^\prime_0+h_1(x_1+x^\prime_1)+\cdots
+h_{n-1}(x_{n-1} +x^\prime_{n-1}) .
\label{2}
\end{equation}
The product of the numbers $u, u^\prime$ is
\begin{equation}
\begin{array}{l}
uu^\prime=x_0 x_0^\prime -x_1x_{n-1}^\prime
-x_2 x_{n-2}^\prime-x_3x_{n-3}^\prime
-\cdots-x_{n-1}x_1^\prime\\
+h_1(x_0 x_1^\prime+x_1x_0^\prime-x_2x_{n-1}^\prime-x_3x_{n-2}^\prime
-\cdots-x_{n-1} x_2^\prime) \\
+h_2(x_0 x_2^\prime+x_1x_1^\prime+x_2x_0^\prime-x_3x_{n-1}^\prime
-\cdots-x_{n-1} x_3^\prime) \\
\vdots\\
+h_{n-1}(x_0 x_{n-1}^\prime+x_1x_{n-2}^\prime
+x_2x_{n-3}^\prime+x_3x_{n-4}^\prime
+\cdots+x_{n-1} x_0^\prime).
\end{array}
\label{3}
\end{equation}
The product $uu^\prime$ can be written as
\begin{equation}
uu^\prime=\sum_{k=0}^{n-1}h_k\sum_{l=0}^{n-1}(-1)^{[(n-k-1+l)/n]}
x_l x^\prime_{k-l+n[(n-k-1+l)/n]}.
\label{3a}
\end{equation}
If $u,u^\prime,u^{\prime\prime}$ are n-complex numbers, the multiplication 
is associative
\begin{equation}
(uu^\prime)u^{\prime\prime}=u(u^\prime u^{\prime\prime})
\label{3b}
\end{equation}
and commutative
\begin{equation}
u u^\prime=u^\prime u ,
\label{3c}
\end{equation}
because the product of the bases, defined in Eq. (\ref{1}), is associative and
commutative. The fact that the multiplication is commutative can be seen also
directly from Eq. (\ref{3}). 
The n-complex zero is $0+h_1\cdot 0+\cdots+h_{n-1}\cdot 0,$ 
denoted simply 0, 
and the n-complex unity is $1+h_1 \cdot 0+\cdots+h_{n-1}\cdot 0,$ 
denoted simply 1.

The inverse of the n-complex number $u=x_0+h_1x_1+h_2x_2+\cdots+h_{n-1}x_{n-1}$
is the n-complex number
$u^\prime
=x^\prime_0+h_1x^\prime_1+h_2x^\prime_2+\cdots+h_{n-1}x^\prime_{n-1}$ 
having the property that
\begin{equation}
uu^\prime=1 .
\label{4}
\end{equation}
Written on components, the condition, Eq. (\ref{4}), is
\begin{equation}
\begin{array}{l}
x_0 x_0^\prime -x_1x_{n-1}^\prime-x_2 x_{n-2}^\prime-x_3x_{n-3}^\prime
-\cdots-x_{n-1}x_1^\prime=1,\\
x_0 x_1^\prime+x_1x_0^\prime-x_2x_{n-1}^\prime-x_3x_{n-2}^\prime
-\cdots-x_{n-1} x_2^\prime=0, \\
x_0 x_2^\prime+x_1x_1^\prime+x_2x_0^\prime-x_3x_{n-1}^\prime
-\cdots-x_{n-1} x_3^\prime=0, \\
\vdots\\
x_0 x_{n-1}^\prime+x_1x_{n-2}^\prime+x_2x_{n-3}^\prime+x_3x_{n-4}^\prime
+\cdots+x_{n-1} x_0^\prime=0.
\end{array}
\label{5}
\end{equation}
The system (\ref{5}) has a solution provided that the determinant of the
system, 
\begin{equation}
\nu={\rm det}(A), 
\label{5b}
\end{equation}
is not equal to zero, $\nu\not=0$, where
\begin{equation}
A=\left(
\begin{array}{ccccc}
x_0     &  -x_{n-1} &   -x_{n-2}  & \cdots  &-x_1\\
x_1     &   x_0     &   -x_{n-1}  & \cdots  &-x_2\\
x_2     &   x_1     &   x_0       & \cdots  &-x_3\\
\vdots  &  \vdots   &  \vdots     & \cdots  &\vdots \\
x_{n-1} &   x_{n-2} &   x_{n-3}   & \cdots  &x_0\\
\end{array}
\right).
\label{6}
\end{equation}
It will be shown that $\nu>0$, and the quantity 
\begin{equation}
\rho=\nu^{1/n}
\label{6b}
\end{equation}
will be called amplitude of the
n-complex number $u=x_0+h_1x_1+h_2x_2+\cdots+h_{n-1}x_{n-1}$.
The quantity $\nu$ can be written as a product of linear factors
\begin{equation}
\nu=\prod_{k=1}^{n}
\left(x_0+\epsilon_k x_1+\epsilon_k^2 x_2+\cdots
+\epsilon^{n-1}_k x_{n-1}\right),
\label{7}
\end{equation}
where $\epsilon_k=e^{i\pi (2k-1)/n}$, $k=1,...,n$, and $i$ being the imaginary
unit. The factors appearing in Eq. (\ref{7}) are of the form 
\begin{equation}
x_0+\epsilon_k x_1+\epsilon_k^2 x_2+\cdots
+\epsilon^{n-1}_k x_{n-1}=v_k+i\tilde v_k,
\label{8}
\end{equation}
where
\begin{equation}
v_k=\sum_{p=0}^{n-1}x_p\cos\frac{\pi (2k-1)p}{n},
\label{9a}
\end{equation}
\begin{equation}
\tilde v_k=\sum_{p=0}^{n-1}x_p\sin\frac{\pi (2k-1)p}{n},
\label{9b}
\end{equation}
for $k=1,...,n$.
The variables $v_k, \tilde v_k, k=1,...,n/2$ will be called canonical polar
n-complex variables. 
It can be seen that $v_k=v_{n-k+1}, \tilde v_k=-\tilde v_{n-k+1}$, for
$k=1,...,n/2$.
Therefore,  the factors appear in Eq. (\ref{7}) in complex-conjugate pairs of
the form $v_k+i\tilde v_k$ and $v_{n-k+1}+i\tilde v_{n-k+1}=v_k-i\tilde v_k$,
where $k=1,...n/2$, so that the determinant $\nu$ is a real and positive
quantity, $\nu>0$, 
\begin{equation}
\nu=\prod_{k=1}^{n/2}\rho_k^2,
\label{9c}
\end{equation}
where
\begin{equation}
\rho_k^2=v_k^2+\tilde v_k^2 .
\label{9d}
\end{equation}
Thus, an n-complex number has an inverse
unless it lies on one of the nodal hypersurfaces 
$\rho_1=0$, or $\rho_2=0$, or ... or $\rho_{n/2}=0$.

\section{Geometric representation of planar n-complex numbers}

The n-complex number $x_0+h_1x_1+h_2x_2+\cdots+h_{n-1}x_{n-1}$
can be represented by 
the point $A$ of coordinates $(x_0,x_1,...,x_{n-1})$. 
If $O$ is the origin of the n-dimensional space,  the
distance from the origin $O$ to the point $A$ of coordinates
$(x_0,x_1,...,x_{n-1})$ has the expression
\begin{equation}
d^2=x_0^2+x_1^2+\cdots+x_{n-1}^2 .
\label{10}
\end{equation}
The quantity $d$ will be called modulus of the n-complex number 
$u=x_0+h_1x_1+h_2x_2+\cdots+h_{n-1}x_{n-1}$. The modulus of an n-complex number
$u$ will be designated by $d=|u|$.

The exponential and trigonometric forms of the n-complex number $u$ can be
obtained conveniently in a rotated system of axes defined by a transformation
which has the form
\begin{eqnarray}
\lefteqn{\left(
\begin{array}{c}
\vdots\\
\xi_k\\
\eta_k\\
\vdots
\end{array}\right)
=\left(
\begin{array}{ccccc}
\vdots&\vdots& &\vdots&\vdots\\
\sqrt{\frac{2}{n}}&\sqrt{\frac{2}{n}}\cos\frac{\pi (2k-1)}{n}&\cdots&
\sqrt{\frac{2}{n}}\cos\frac{\pi (2k-1)(n-2)}{n}&\sqrt{\frac{2}{n}}
\cos\frac{\pi (2k-1)(n-1)}{n}\\
0&\sqrt{\frac{2}{n}}\sin\frac{\pi (2k-1)}{n}&
\cdots&\sqrt{\frac{2}{n}}\sin\frac{\pi (2k-1)(n-2)}{n}&
\sqrt{\frac{2}{n}}\sin\frac{\pi (2k-1)(n-1)}{n}\\
\vdots&\vdots&&\vdots&\vdots
\end{array}
\right)
\left(
\begin{array}{c}
x_0\\
\vdots\\ 
\vdots\\
x_{n-1}
\end{array}
\right),\nonumber}\\
&&
\label{11}
\end{eqnarray}
where $k=1, 2, ... , n/2$.
The lines of the matrices in Eq. (\ref{11}) give the components
of the $n$ vectors of the new basis system of axes. These vectors have unit
length and are orthogonal to each other.
By comparing Eqs. (\ref{9a})-(\ref{9b}) and (\ref{11}) it can be
seen that
\begin{equation}
v_k= \sqrt{\frac{n}{2}}\xi_k , \tilde v_k= \sqrt{\frac{n}{2}}\eta_k ,
\label{12b}
\end{equation}
i.e. the two sets of variables differ only by a scale factor.

The sum of the squares of the variables $v_k,\tilde v_k$ is
\begin{equation}
\sum_{k=1}^{n/2}(v_k^2+\tilde v_k^2)=\frac{n}{2}d^2.
\label{13}
\end{equation}
The relation (\ref{13}) has been obtained with the aid of the relation
\begin{equation}
\sum_{k=1}^{n/2}\cos\frac{\pi (2k-1)p}{n}=0, 
\label{15}
\end{equation}
for $p=1,...,n-1$.
From Eq. (\ref{13}) it results that
\begin{equation}
d^2=\frac{2}{n}\sum_{k=1}^{n/2}\rho_k^2 .
\label{17}
\end{equation}
The relation (\ref{17}) shows that the square of the distance
$d$, Eq. (\ref{10}), is equal to the sum of the squares of the projections
$\rho_k\sqrt{2/n}$. This is consistent with the fact that the transformation in
Eq. (\ref{11}) is unitary.

The position of the point $A$ of coordinates $(x_0,x_1,...,x_{n-1})$ can be
also described with the aid of the distance $d$, Eq. (\ref{10}), and of $n-1$
angles defined further. Thus, in the plane of the axes $v_k,\tilde v_k$, the
azimuthal angle $\phi_k$ can be introduced by the relations 
\begin{equation}
\cos\phi_k=v_k/\rho_k,\:\sin\phi_k=\tilde v_k/\rho_k, 
\label{19a}
\end{equation}
where $0\leq \phi_k<2\pi, \;k=1,...,n/2$, so that there are $n/2$ azimuthal
angles. 
The radial distance $\rho_k$ in the plane of the axes $v_k,\tilde v_k$ has been
defined in Eq. (\ref{9d}).
If the projection of the point $A$ on the plane of the axes $v_k,\tilde v_k$ is
$A_k$, and the projection of the point $A$ on the 4-dimensional space defined
by the axes $v_1, \tilde v_1, v_k,\tilde v_k$ is $A_{1k}$, the angle
$\psi_{k-1}$ between the line $OA_{1k}$ and the 2-dimensional plane defined by
the axes $v_k,\tilde v_k$ is  
\begin{equation}
\tan\psi_{k-1}=\rho_1/\rho_k, 
\label{19b}
\end{equation}
where $0\leq\psi_k\leq\pi/2, k=2,...,n/2$,
so that there are $n/2-1$ planar angles. Thus, the
position of the point $A$ is described by the distance $d$, by $n/2$
azimuthal angles and by $n/2-1$ 
planar angles. These angles are shown in Fig. 2. 

The variables $\rho_k$ can be expressed in terms of $d$ and the planar angles
$\psi_k$ as
\begin{equation}
\rho_k=\frac{\rho_1}{\tan\psi_{k-1}}, 
\label{20a}
\end{equation}
for $k=2,...,n/2$, where
\begin{eqnarray}
\rho_1^2=\frac{nd^2}{2}
\left(1
+\frac{1}{\tan^2\psi_1}+\frac{1}{\tan^2\psi_2}+\cdots
+\frac{1}{\tan^2\psi_{n/2-1}}\right)^{-1}.
\label{20b}
\end{eqnarray}

If
$u^\prime=x_0^\prime+h_1x_1^\prime+h_2x_2^\prime+\cdots+h_{n-1}x_{n-1}^\prime, 
u^{\prime\prime}=x^{\prime\prime}_0+h_1x^{\prime\prime}_1
+h_2x^{\prime\prime}_2+\cdots+h_{n-1}x^{\prime\prime}_{n-1}$ 
are n-complex numbers of parameters
$\rho_k^\prime,\psi_k^\prime,\phi_k^\prime$ and respectively  
$\rho_k^{\prime\prime}, \psi_k^{\prime\prime},
\phi_k^{\prime\prime}$, then the parameters
$v_+,\rho_k,\psi_k,\phi_k$ of the product n-complex number
$u=u^\prime u^{\prime\prime}$ are given by
\begin{equation}
\rho_k=\rho_k^\prime\rho_k^{\prime\prime}, 
\label{21b}
\end{equation}
for $k=1,..., n/2$,
\begin{equation}
\tan\psi_k=\tan\psi_k^\prime \tan\psi_k^{\prime\prime}, 
\label{21d}
\end{equation}
for $k=1,...,n/2-1$, 
\begin{equation}
\phi_k=\phi_k^\prime+\phi_k^{\prime\prime}, 
\label{21e}
\end{equation}
for $k=1,..., n/2$.
The Eqs. (\ref{21b})-(\ref{21e}) are a consequence of the relations
\begin{equation}
v_k=v_k^\prime v_k^{\prime\prime}-\tilde v_k^\prime \tilde
v_k^{\prime\prime},\; 
\tilde v_k=v_k^\prime \tilde v_k^{\prime\prime}+\tilde v_k^\prime
v_k^{\prime\prime}, 
\label{22}
\end{equation}
and of the corresponding relations of definition. Then the product $\nu$ in
Eq. (\ref{9c}) has the property that
\begin{equation}
\nu=\nu^\prime\nu^{\prime\prime} ,
\label{23}
\end{equation}
and the amplitude $\rho$ defined in Eq. (\ref{6b}) has the property that
\begin{equation}
\rho=\rho^\prime\rho^{\prime\prime} .
\label{24}
\end{equation}

The fact that the amplitude of the product is equal to the product of the 
amplitudes, as written in Eq. (\ref{24}), can 
be demonstrated also by using a representation of the 
n-complex numbers by matrices, in which the n-complex number 
$u=x_0+h_1x_1+h_2x_2+\cdots+h_{n-1}x_{n-1}$ is represented by the matrix
\begin{equation}
U=\left(
\begin{array}{ccccc}
x_0      &   x_1      &   x_2    & \cdots  &x_{n-1}\\
-x_{n-1} &   x_0      &   x_1    & \cdots  &x_{n-2}\\
-x_{n-2} &  - x_{n-1} &   x_0    & \cdots  &x_{n-3}\\
\vdots   &  \vdots    &  \vdots  & \cdots  &\vdots \\
-x_1     &   -x_2     &  - x_3   & \cdots  &x_0\\
\end{array}
\right).
\label{24a}
\end{equation}
The product $u=u^\prime u^{\prime\prime}$ is be
represented by the matrix multiplication $U=U^\prime U^{\prime\prime}$.
The relation (\ref{23}) is then a consequence of the fact the determinant 
of the product of matrices is equal to the product of the determinants 
of the factor matrices. 
The use of the representation  of the n-complex numbers with matrices
provides an alternative demonstration of the fact that the product of
n-complex numbers is associative, as stated in 
Eq. (\ref{3b}).

According to Eqs. (\ref{13} and (\ref{9d}), the modulus of the product
$uu^\prime$ is given by  
\begin{equation}
|uu^\prime|^2=
\frac{2}{n}\sum_{k=1}^{n/2}(\rho_k\rho_k^\prime)^2 .
\label{25a}
\end{equation}
Thus, if the product of two n-complex numbers is zero, $uu^\prime=0$, then
$\rho_k\rho_k^\prime=0, k=1,...,n/2$. This means that either $u=0$, or
$u^\prime=0$, or $u, u^\prime$ belong to orthogonal hypersurfaces in such a way
that the afore-mentioned products of components should be equal to zero.

\section{The planar n-dimensional cosexponential functions}

The exponential function of the n-complex variable $u$ can be defined by the
series
\begin{equation}
\exp u = 1+u+u^2/2!+u^3/3!+\cdots . 
\label{26}
\end{equation}
It can be checked by direct multiplication of the series that
\begin{equation}
\exp(u+u^\prime)=\exp u \cdot \exp u^\prime . 
\label{27}
\end{equation}
If $u=x_0+h_1x_1+h_2x_2+\cdots+h_{n-1}x_{n-1}$,
then $\exp u$ can be calculated as 
$\exp u=\exp x_0 \cdot \exp (h_1x_1) \cdots \exp (h_{n-1}x_{n-1})$.

It can be seen with the aid of the representation in Fig. 1 that 
\begin{equation}
h_k^{n+p}=(-1)^k h_k^p, \:p\:\:{\rm integer},
\label{28}
\end{equation}
where $k=1,...,n-1$. For $k$ even, $e^{h_k y}$ can be written as
\begin{equation}
e^{h_k y}=\sum_{p=0}^{n-1}(-1)^{[kp/n]}h_{kp-n[kp/n]}g_{np}(y),
\label{28bx}
\end{equation}
where $h_0=1$, and where $g_{np}$ are the polar n-dimensional cosexponential
functions. \cite{1a} For odd $k$,  $e^{h_k y}$ is
\begin{equation}
e^{h_k y}=\sum_{p=0}^{n-1}(-1)^{[kp/n]}h_{kp-n[kp/n]}f_{np}(y),
\label{28b}
\end{equation}
where the functions $f_{nk}$, which will be
called planar cosexponential functions in $n$ dimensions, are
\begin{equation}
f_{nk}(y)=\sum_{p=0}^\infty (-1)^p \frac{y^{k+pn}}{(k+pn)!}, 
\label{29}
\end{equation}
for $ k=0,1,...,n-1$.

The planar cosexponential functions  of even index $k$ are
even functions, $f_{n,2l}(-y)=f_{n,2l}(y)$,  
and the planar cosexponential functions of odd index 
are odd functions, $f_{n,2l+1}(-y)=-f_{n,2l+1}(y)$, $l=0,...,n/2-1$ . 

The planar n-dimensional cosexponential function $f_{nk}(y)$ is related to the
polar n-dimensional cosexponential function $g_{nk}(y)$ discussed in \cite{1a}
by 
the relation 
\begin{equation}
f_{nk}(y)=e^{-i\pi k/n}g_{nk}\left(e^{i\pi/n}y\right), 
\label{30a}
\end{equation}
for $k=0,...,n-1$.
The expression of the planar n-dimensional cosexponential functions is then
\begin{equation}
f_{nk}(y)=\frac{1}{n}\sum_{l=1}^{n}
\exp\left[y\cos\left(\frac{\pi (2l-1)}{n}\right)
\right]
\cos\left[y\sin\left(\frac{\pi (2l-1)}{n}\right)-\frac{\pi (2l-1)k}{n}\right], 
\label{30}
\end{equation}
for $k=0,1,...,n-1$.
The planar cosexponential function defined in Eq. (\ref{29}) has the expression
given in Eq. (\ref{30}) for any natural value of $n$, this result not being
restricted to even values of $n$.
In order to check that the function in Eq. (\ref{30}) has the series expansion
written in Eq. (\ref{29}), the right-hand side of Eq. (\ref{30}) will be
written as
\begin{equation}
f_{nk}(y)=\frac{1}{n}\sum_{l=1}^{n}{\rm Re}\left\{
\exp\left[\left(\cos\frac{\pi (2l-1)}{n}+i\sin\frac{\pi (2l-1)}{n}\right)y
-i\frac{\pi k(2l-1)}{n}\right]\right\}, 
\label{31}
\end{equation}
for $k=0,1,...,n-1$, where ${\rm Re} (a+ib)=a$, with $a$ and $b$ real numbers.
The part of the 
exponential depending on $y$ can be expanded in a series,
\begin{equation}
f_{nk}(y)=\frac{1}{n}\sum_{p=0}^\infty\sum_{l=1}^{n}{\rm Re}
\left\{\frac{1}{p!}
\exp\left[i\frac{\pi (2l-1)}{n}(p-k)\right]y^p\right\}, 
\label{32}
\end{equation}
for $k=0,1,...,n-1$.
The expression of $f_{nk}(y)$ becomes
\begin{equation}
f_{nk}(y)=\frac{1}{n}\sum_{p=0}^\infty\sum_{l=1}^{n}
\left\{\frac{1}{p!}
\cos\left[\frac{\pi (2l-1)}{n}(p-k)\right]y^p\right\}, 
\label{33}
\end{equation}
where $k=0,1,...,n-1$ and, since
\begin{equation}
\frac{1}{n}\sum_{l=1}^{n}\cos\left[\frac{\pi (2l-1)}{n}(p-k)\right]
=\left\{
\begin{array}{l}
1, \:\:{\rm if}\:\: p-k\:\: {\rm is \:\;an\:\;even\:\;multiple\:\;of\:\;}n,\\
-1, \:\:{\rm if}\:\: p-k\:\: {\rm is \:\;an\:\;odd\:\;multiple\:\;of\:\;}n,\\
0, \:\:{\rm otherwise},
\end{array}
\right.
\label{34}
\end{equation}
this yields indeed the expansion in Eq. (\ref{29}).

It can be shown from Eq. (\ref{30}) that
\begin{equation}
\sum_{k=0}^{n-1}f_{nk}^2(y)=\frac{1}{n}\sum_{l=1}^{n}\exp\left[2y\cos\left(
\frac{\pi (2l-1)}{n}\right)\right].
\label{34a}
\end{equation}
It can be seen that the right-hand side of Eq. (\ref{34a}) does not contain
oscillatory terms. If $n$ is a multiple of 4, it can be shown by replacing $y$
by $iy$ in Eq. (\ref{34a}) that
\begin{equation}
\sum_{k=0}^{n-1}(-1)^kf_{nk}^2(y)=\frac{4}{n}
\sum_{l=1}^{n/4}\cos\left[2y\cos\left(\frac{\pi (2l-1)}{n}\right)\right],
\label{34b}
\end{equation}
which does not contain exponential terms.

For odd $n$, the planar n-dimensional cosexponential function $f_{nk}(y)$ is
related to the n-dimensional cosexponential function $g_{nk}(y)$ discussed in
\cite{1a} also by the relation 
\begin{equation}
f_{nk}(y)=(-1)^kg_{nk}(-y),
\label{34c}
\end{equation}
as can be seen by comparing the series for the two classes of functions.
For values of the form $n=4p+2$, $p=0,1,2,...$, the planar n-dimensional
cosexponential function $f_{nk}(y)$ is related to the n-dimensional
cosexponential function $g_{nk}(y)$ by the relation  
\begin{equation}
f_{nk}(y)=e^{-i\pi k/2}g_{nk}(iy).
\label{34cx}
\end{equation}

Addition theorems for the planar n-dimensional cosexponential functions can be
obtained from the relation $\exp h_1(y+z)=\exp h_1 y \cdot\exp h_1 z $, by
substituting the expression of the exponentials as given in Eq. (\ref{28b})
for $k=1$, $e^{h_1 y}=f_{n0}(y)+h_1f_{n1}(y)+\cdots+h_{n-1} f_{n,n-1}(y)$,
\begin{eqnarray}
\lefteqn{f_{nk}(y+z)=f_{n0}(y)f_{nk}(z)
+f_{n1}(y)f_{n,k-1}(z)+\cdots+f_{nk}(y)f_{n0}(z)\nonumber}\\
&&-f_{n,k+1}(y)f_{n,n-1}(z)-f_{n,k+2}(y)f_{n,n-2}(z)
-\cdots-f_{n,n-1}(y)f_{n,k+1}(z) ,
\nonumber\\
&&
\label{35a}
\end{eqnarray}
where $k=0,1,...,n-1$.
For $y=z$ the relations (\ref{35a}) take the form
\begin{eqnarray}
\lefteqn{f_{nk}(2y)=f_{n0}(y)f_{nk}(y)+f_{n1}(y)f_{n,k-1}(y)
+\cdots+f_{nk}(y)f_{n0}(y)\nonumber}\\
&&-f_{n,k+1}(y)f_{n,n-1}(y)-f_{n,k+2}(y)f_{n,n-2}(y)
-\cdots-f_{n,n-1}(y)f_{n,k+1}(y) ,
\nonumber\\
&&
\label{35b}
\end{eqnarray}
where $k=0,1,...,n-1$.
For $y=-z$ the relations (\ref{35a}) and (\ref{29}) yield
\begin{equation}
f_{n0}(y)f_{n0}(-y)
-f_{n1}(y)f_{n,n-1}(-y)-f_{n2}(y)f_{n,n-2}(-y)
-\cdots-f_{n,n-1}(y)f_{n1}(-y)=1 ,
\label{36a}
\end{equation}
\begin{eqnarray}
\lefteqn{f_{n0}(y)f_{nk}(-y)+f_{n1}(y)f_{n,k-1}(-y)
+\cdots+f_{nk}(y)f_{n0}(-y)\nonumber}\\
&&-f_{n,k+1}(y)f_{n,n-1}(-y)-f_{n,k+2}(y)f_{n,n-2}(-y)
-\cdots-f_{n,n-1}(y)f_{n,k+1}(-y)=0 ,
\nonumber\\
&&
\label{36b}
\end{eqnarray}
where $k=1,...,n-1$.

From Eq. (\ref{28bx}) it can be shown, for even $k$ and 
natural numbers $l$, that
\begin{equation}
\left(\sum_{p=0}^{n-1}(-1)^{[kp/n]}h_{kp-n[kp/n]}g_{np}(y)\right)^l
=\sum_{p=0}^{n-1}(-1)^{[kp/n]}h_{kp-n[kp/n]}g_{np}(ly), 
\label{37x}
\end{equation}
where $k=0,1,...,n-1$.
For odd $k$ and natural numbers $l$, Eq. (\ref{28b}) implies
\begin{equation}
\left(\sum_{p=0}^{n-1}(-1)^{[kp/n]}h_{kp-n[kp/n]}f_{np}(y)\right)^l
=\sum_{p=0}^{n-1}(-1)^{[kp/n]}h_{kp-n[kp/n]}f_{np}(ly), 
\label{37}
\end{equation}
where $k=0,1,...,n-1$.
For $k=1$ the relation (\ref{37}) is
\begin{equation}
\left\{f_{n0}(y)+h_1f_{n1}(y)+\cdots+h_{n-1}f_{n,n-1}(y)\right\}^l
=f_{n0}(ly)+h_1f_{n1}(ly)+\cdots+h_{n-1}f_{n,n-1}(ly).
\label{37b}
\end{equation}

If
\begin{equation}
a_k=\sum_{p=0}^{n-1}f_{np}(y)\cos\left(\frac{\pi (2k-1)p}{n}\right), 
\label{38a}
\end{equation}
and
\begin{equation}
b_k=\sum_{p=0}^{n-1}f_{np}(y)\sin\left(\frac{\pi (2k-1)p}{n}\right), 
\label{38b}
\end{equation}
for $k=1,...,n$, where $f_{np}(y)$ are the planar cosexponential functions in
Eq. (\ref{30}), it can be shown that
\begin{equation}
a_k=\exp\left[y\cos\left(\frac{\pi (2k-1)}{n}\right)\right]
\cos\left[y\sin\left(\frac{\pi (2k-1)}{n}\right)\right], 
\label{39a}
\end{equation}
\begin{equation}
b_k=\exp\left[y\cos\left(\frac{\pi (2k-1)}{n}\right)\right]
\sin\left[y\sin\left(\frac{\pi (2k-1)}{n}\right)\right], 
\label{39b}
\end{equation}
for $k=1,...,n$. If
\begin{equation}
G_k^2=a_k^2+b_k^2, 
\label{40}
\end{equation}
from Eqs. (\ref{39a}) and (\ref{39b}) it results that
\begin{equation}
G_k^2=\exp\left[2y\cos\left(\frac{\pi (2k-1)}{n}\right)\right], 
\label{41}
\end{equation}
for $k=1,...,n$.
Then the planar n-dimensional cosexponential functions have the property that
\begin{equation}
\prod_{p=1}^{n/2}G_p^2=1.
\label{43a}
\end{equation}

The planar n-dimensional cosexponential functions are solutions of the
$n^{\rm th}$-order differential equation
\begin{equation}
\frac{d^n\zeta}{du^n}=-\zeta ,
\label{44}
\end{equation}
whose solutions are of the form
$\zeta(u)=A_0f_{n0}(u)+A_1f_{n1}(u)+\cdots+A_{n-1}f_{n,n-1}(u).$ 
It can be checked that the derivatives of the planar cosexponential functions
are related by
\begin{equation}
\frac{df_{n0}}{du}=-f_{n,n-1}, \:
\frac{df_{n1}}{du}=f_{n0}, \:...,
\frac{df_{n,n-2}}{du}=f_{n,n-3} ,
\frac{df_{n,n-1}}{du}=f_{n,n-2} .
\label{45}
\end{equation}

\section{Exponential and trigonometric forms of planar n-complex numbers}

In order to obtain the exponential and trigonometric forms of n-complex
numbers, a canonical base
$e_1,\tilde e_1,...,e_{n/2},\tilde e_{n/2}$ for the planar n-complex numbers
will be introduced by the relations 
\begin{equation}
\left(
\begin{array}{c}
\vdots\\
e_k\\
\tilde e_k\\
\vdots
\end{array}\right)
=\left(
\begin{array}{ccccc}
\vdots&\vdots& &\vdots&\vdots\\
\frac{2}{n}&\frac{2}{n}\cos\frac{\pi (2k-1)}{n}&
\cdots&\frac{2}{n}\cos\frac{\pi (2k-1)(n-2)}{n}&
\frac{2}{n}\cos\frac{\pi (2k-1)(n-1)}{n}\\
0&\frac{2}{n}\sin\frac{\pi (2k-1)}{n}&\cdots&
\frac{2}{n}\sin\frac{\pi (2k-1)(n-2)}{n}&\frac{2}{n}
\sin\frac{\pi (2k-1)(n-1)}{n}\\
\vdots&\vdots&&\vdots&\vdots
\end{array}
\right)
\left(
\begin{array}{c}
1\\h_1\\
\vdots\\
h_{n-1}
\end{array}
\right),
\label{e11}
\end{equation}
where $k=1, 2, ... , n/2$.

The multiplication relations for the bases $e_k, \tilde e_k$ are
\begin{eqnarray}
e_k^2=e_k, \tilde e_k^2=-e_k, e_k \tilde e_k=\tilde e_k , e_ke_l=0, e_k\tilde
e_l=0, \tilde e_k\tilde e_l=0, k\not=l, 
\label{e12a}
\end{eqnarray}
where $k,l=1,...,n/2$.
The moduli of the bases $e_k, \tilde e_k$ are
\begin{equation}
|e_k|=\sqrt{\frac{2}{n}}, |\tilde e_k|=\sqrt{\frac{2}{n}}.
\label{e12c}
\end{equation}
It can be shown that
\begin{eqnarray}
x_0+h_1x_1+\cdots+h_{n-1}x_{n-1}= 
\sum_{k=1}^{n/2}(e_k v_k+\tilde e_k \tilde v_k).
\label{e13a}
\end{eqnarray}
The relation (\ref{e13a} gives the canonical form of a planar n-complex
number.

Using the properties of the bases in Eqs. (\ref{e11}) it can
be shown that
\begin{equation}
\exp(\tilde e_k\phi_k)=1-e_k+e_k\cos\phi_k+\tilde e_k\sin\phi_k ,
\label{46a}
\end{equation}
\begin{equation}
\exp(e_k\ln\rho_k)=1-e_k+e_k\rho_k ,
\label{46b}
\end{equation}
By multiplying the relations (\ref{46a}), (\ref{46b}) it results that 
\begin{equation}
\exp\left[\sum_{k=1}^{n/2}
(e_k\ln \rho_k+\tilde e_k\phi_k)\right] 
=\sum_{k=1}^{n/2}(e_k v_k+\tilde e_k \tilde v_k),
\label{47a}
\end{equation}
where the fact has ben used that 
\begin{equation}
\sum_{k=1}^{n/2} e_k =1, 
\label{47b}
\end{equation}
the latter relation being a consequence of Eqs. (\ref{e11}) and (\ref{15}).

By comparing Eqs. (\ref{e13a}) and (\ref{47a}), it can be seen that
\begin{eqnarray}
x_0+h_1x_1+\cdots+h_{n-1}x_{n-1}=
\exp\left[\sum_{k=1}^{n/2}
(e_k\ln \rho_k+\tilde e_k\phi_k)\right] .
\label{49a}
\end{eqnarray}
Using the expression of the bases in Eq. (\ref{e11}) yields
the exponential form of the n-complex number
$u=x_0+h_1x_1+\cdots+h_{n-1}x_{n-1}$ as
\begin{eqnarray}
u=\rho\exp\left\{\sum_{p=1}^{n-1}h_p\left[
-\frac{2}{n}\sum_{k=2}^{n/2}
\cos\left(\frac{\pi (2k-1)p}{n}\right)\ln\tan\psi_{k-1}
\right]
+\sum_{k=1}^{n/2}\tilde e_k\phi_k
\right\},
\label{50a}
\end{eqnarray}
where $\rho$ is the amplitude defined in Eq. (\ref{6b}), and has according to
Eq. (\ref{9c}) the expression
\begin{equation}
\rho=\left(\rho_1^2\cdots\rho_{n/2}^2\right)^{1/n}.
\label{50aa}
\end{equation}

It can be checked with the aid of Eq. (\ref{46a}) that
the n-complex number $u$ can also be written as
\begin{eqnarray}
x_0+h_1x_1+\cdots+h_{n-1}x_{n-1}
=\left(\sum_{k=1}^{n/2}e_k \rho_k\right)
\exp\left(\sum_{k=1}^{n/2}\tilde e_k\phi_k\right).
\label{51a}
\end{eqnarray}
Writing in Eq. (\ref{51a}) the radius $\rho_1$, Eq.
(\ref{20b}), as a factor and expressing the 
variables in terms of the planar
angles with the aid of Eq. (\ref{19b}) yields the
trigonometric form of the n-complex number $u$ as
\begin{eqnarray}
\lefteqn{u=d
\left(\frac{n}{2}\right)^{1/2}
\left(1+\frac{1}{\tan^2\psi_1}+\frac{1}{\tan^2\psi_2}+\cdots
+\frac{1}{\tan^2\psi_{n/2-1}}\right)^{-1/2}\nonumber}\\
&&\left(e_1+\sum_{k=2}^{n/2}\frac{e_k}{\tan\psi_{k-1}}\right)
\exp\left(\sum_{k=1}^{n/2}\tilde e_k\phi_k\right).
\label{52a}
\end{eqnarray}
In Eq. (\ref{52a}),  the n-complex
number $u$, written in trigonometric form, is the
product of the modulus $d$, of a part depending on the planar
angles $\psi_1,...,\psi_{n/2-1}$, and of a factor depending
on the azimuthal angles $\phi_1,...,\phi_{n/2}$. 
Although the modulus of a product of n-complex numbers is not equal in
general to the product of the moduli of the factors,
it can be checked that the modulus of the factors in Eq. (\ref{52a}) are
\begin{eqnarray}
\lefteqn{
\left|e_1+\sum_{k=2}^{n/2}\frac{e_k}{\tan\psi_{k-1}}\right|\nonumber}\\
&&=\left(\frac{2}{n}\right)^{1/2}
\left(1+\frac{1}{\tan^2\psi_1}+\frac{1}{\tan^2\psi_2}+\cdots
+\frac{1}{\tan^2\psi_{n/2-1}}\right)^{1/2},
\label{52c}
\end{eqnarray}
and
\begin{eqnarray}
\left|\exp\left(\sum_{k=1}^{n/2}\tilde e_k\phi_k\right)\right|=1.
\label{52e}
\end{eqnarray}

The modulus $d$ in Eqs. (\ref{52a}) can be expressed in terms
of the amplitude $\rho$ as
\begin{eqnarray}
\lefteqn{d=\rho \frac{2^{(n-2)/2n}}{\sqrt{n}}
\left(\tan\psi_1\cdots\tan\psi_{n/2-1}\right)^{2/n}\nonumber}\\
&&\left(1+\frac{1}{\tan^2\psi_1}+\frac{1}{\tan^2\psi_2}+\cdots
+\frac{1}{\tan^2\psi_{n/2-1}}\right)^{1/2}.
\label{53a}
\end{eqnarray}

\section{Elementary functions of a planar n-complex variable}

The logarithm $u_1$ of the n-complex number $u$, $u_1=\ln u$, can be defined
as the solution of the equation
\begin{equation}
u=e^{u_1} .
\label{54}
\end{equation}
The relation (\ref{47a}) shows that $\ln u$ exists  
as an n-complex function with real components for all values of
$x_0,...,x_{n-1}$ for which $\rho\not=0$.
The expression of the logarithm, obtained from Eq. (\ref{49a}) is
\begin{equation}
\ln u=\sum_{k=1}^{n/2}
(e_k\ln \rho_k+\tilde e_k\phi_k).
\label{55a}
\end{equation}
An expression of the logarithm depending on the amplitude $\rho$ can be
obtained from the exponential forms in Eq. (\ref{50a}) as
\begin{eqnarray}
\ln u=\ln \rho+\sum_{p=1}^{n-1}h_p\left[
-\frac{2}{n}\sum_{k=2}^{n/2}
\cos\left(\frac{\pi (2k-1)p}{n}\right)\ln\tan\psi_{k-1}
\right]+\sum_{k=1}^{n/2}\tilde e_k\phi_k.
\label{56a}
\end{eqnarray}

The function $\ln u$ is multivalued because of the presence of the terms 
$\tilde e_k\phi_k$.
It can be inferred from Eqs. (\ref{21b})-(\ref{21e}) and (\ref{24}) that
\begin{equation}
\ln(uu^\prime)=\ln u+\ln u^\prime ,
\label{57}
\end{equation}
up to integer multiples of $2\pi\tilde e_k, k=1,...,n/2$.

The power function $u^m$ can be defined for real values of $m$ as
\begin{equation}
u^m=e^{m\ln u} .
\label{58}
\end{equation}
Using the expression of $\ln u$ in Eq. (\ref{55a}) yields
\begin{equation}
u^m=\sum_{k=1}^{n/2}
\rho_k^m(e_k\cos m\phi_k+\tilde e_k\sin m\phi_k).
\label{59a}
\end{equation}
The power function is multivalued unless $m$ is an integer. 
For integer $m$, it can be inferred from Eq. (\ref{57}) that
\begin{equation}
(uu^\prime)^m=u^m\:u^{\prime m} .
\label{59}
\end{equation}

The trigonometric functions $\cos u$ and $\sin u$ of an n-complex variable $u$
are defined by the series
\begin{equation}
\cos u = 1 - u^2/2!+u^4/4!+\cdots, 
\label{60}
\end{equation}
\begin{equation}
\sin u=u-u^3/3!+u^5/5! +\cdots .
\label{61}
\end{equation}
It can be checked by series multiplication that the usual addition theorems
hold for the n-complex numbers $u, u^\prime$,
\begin{equation}
\cos(u+u^\prime)=\cos u\cos u^\prime - \sin u \sin u^\prime ,
\label{62}
\end{equation}
\begin{equation}
\sin(u+u^\prime)=\sin u\cos u^\prime + \cos u \sin u^\prime .
\label{63}
\end{equation}

In order to obtain expressions for the trigonometric functions of n-complex
variables, these will be expressed with the aid of the imaginary unit $i$ as
\begin{equation}
\cos u=\frac{1}{2}(e^{iu}+e^{-iu}),\:\sin u=\frac{1}{2i}(e^{iu}-e^{-iu}).
\label{64}
\end{equation}
The imaginary unit $i$ is used for the convenience of notations, and it does
not appear in the final results.
The validity of Eq. (\ref{64}) can be checked by comparing the series for the
two sides of the relations.
Since the expression of the exponential function $e^{h_k y}$ in terms of the
units $1, h_1, ... h_{n-1}$ given in Eq. (\ref{28b}) depends on the planar
cosexponential functions $f_{np}(y)$, the expression of the trigonometric
functions will depend on the functions
$f_{p+}^{(c)}(y)=(1/2)[f_{np}(iy)+f_{np}(-iy)]$ and 
$f_{p-}^{(c)}(y)=(1/2i)[f_{np}(iy)-f_{np}(-iy)]$, 
\begin{equation}
\cos(h_k y)=\sum_{p=0}^{n-1}(-1)^{[kp/n]}h_{kp-n[kp/n]}f_{p+}^{(c)}(y),
\label{66a}
\end{equation}
\begin{equation}
\sin(h_k y)=\sum_{p=0}^{n-1}(-1)^{[kp/n]}h_{kp-n[kp/n]}f_{p-}^{(c)}(y),
\label{66b}
\end{equation}
where
\begin{eqnarray}
\lefteqn{f_{p+}^{(c)}(y)=\frac{1}{n}\sum_{l=1}^{n}\left\{
\cos\left[y\cos\left(\frac{\pi (2l-1)}{n}\right)\right]
\cosh\left[y\sin\left(\frac{\pi (2l-1)}{n}\right)\right]
\cos\left(\frac{\pi (2l-1)p}{n}\right)\right.\nonumber}\\
&&\left.-\sin\left[y\cos\left(\frac{\pi (2l-1)}{n}\right)\right]
\sinh\left[y\sin\left(\frac{\pi (2l-1)}{n}\right)\right]
\sin\left(\frac{\pi (2l-1)p}{n}\right)
\right\},
\label{65a}
\end{eqnarray}
\begin{eqnarray}
\lefteqn{f_{p-}^{(c)}(y)=\frac{1}{n}\sum_{l=1}^{n}\left\{
\sin\left[y\cos\left(\frac{\pi (2l-1)}{n}\right)\right]
\cosh\left[y\sin\left(\frac{\pi (2l-1)}{n}\right)\right]
\cos\left(\frac{\pi (2l-1)p}{n}\right)\right.\nonumber}\\
&&\left.+\cos\left[y\cos\left(\frac{\pi (2l-1)}{n}\right)\right]
\sinh\left[y\sin\left(\frac{\pi (2l-1)}{n}\right)\right]
\sin\left(\frac{\pi (2l-1)p}{n}\right)
\right\}.
\label{65b}
\end{eqnarray}

The hyperbolic functions $\cosh u$ and $\sinh u $ of the n-complex variable
$u$ can be defined by the series
\begin{equation}
\cosh u = 1 + u^2/2!+u^4/4!+\cdots, 
\label{66}
\end{equation}
\begin{equation}
\sinh u=u+u^3/3!+u^5/5! +\cdots .
\label{67}
\end{equation}
It can be checked by series multiplication that the usual addition theorems
hold for the n-complex numbers $u, u^\prime$,
\begin{equation}
\cosh(u+u^\prime)=\cosh u\cosh u^\prime + \sinh u \sinh u^\prime ,
\label{68}
\end{equation}
\begin{equation}
\sinh(u+u^\prime)=\sinh u\cosh u^\prime + \cosh u \sinh u^\prime .
\label{69}
\end{equation}
In order to obtain expressions for the hyperbolic functions of n-complex
variables, these will be expressed as
\begin{equation}
\cosh u=\frac{1}{2}(e^{u}+e^{-u}),\:\sinh u=\frac{1}{2}(e^{u}-e^{-u}).
\label{70}
\end{equation}
The validity of Eq. (\ref{70}) can be checked by comparing the series for the
two sides of the relations.
Since the expression of the exponential function $e^{h_k y}$ in terms of the
units $1, h_1, ... h_{n-1}$ given in Eq. (\ref{28b}) depends on the planar
cosexponential functions $f_{np}(y)$, the expression of the hyperbolic
functions will depend on the even part $f_{p+}(y)=(1/2)[f_{np}(y)+f_{np}(-y)]$
and on 
the odd part $f_{p-}(y)=(1/2)[f_{np}(y)-f_{np}(-y)]$ of $f_{np}$, 
\begin{equation}
\cosh(h_k y)=\sum_{p=0}^{n-1}(-1)^{[kp/n]}h_{kp-n[kp/n]}f_{p+}(y),
\label{71a}
\end{equation}
\begin{equation}
\sinh(h_k y)=\sum_{p=0}^{n-1}(-1)^{[kp/n]}h_{kp-n[kp/n]}f_{p-}(y),
\label{71b}
\end{equation}
where
\begin{eqnarray}
\lefteqn{f_{p+}(y)=\frac{1}{n}\sum_{(2l-1)=1}^{n}\left\{
\cosh\left[y\cos\left(\frac{\pi (2l-1)}{n}\right)\right]
\cos\left[y\sin\left(\frac{\pi (2l-1)}{n}\right)\right]
\cos\left(\frac{\pi (2l-1)p}{n}\right)\right.\nonumber}\\
&&\left.+\sinh\left[y\cos\left(\frac{\pi (2l-1)}{n}\right)\right]
\sin\left[y\sin\left(\frac{\pi (2l-1)}{n}\right)\right]
\sin\left(\frac{\pi (2l-1)p}{n}\right)
\right\},
\label{72a}
\end{eqnarray}
\begin{eqnarray}
\lefteqn{f_{p-}(y)=\frac{1}{n}\sum_{l=1}^{n}\left\{
\sinh\left[y\cos\left(\frac{\pi (2l-1)}{n}\right)\right]
\cos\left[y\sin\left(\frac{\pi (2l-1)}{n}\right)\right]
\cos\left(\frac{\pi (2l-1)p}{n}\right)\right.\nonumber}\\
&&\left.+\cosh\left[y\cos\left(\frac{\pi (2l-1)}{n}\right)\right]
\sin\left[y\sin\left(\frac{\pi (2l-1)}{n}\right)\right]
\sin\left(\frac{\pi (2l-1)p}{n}\right)
\right\}.
\label{72b}
\end{eqnarray}

The exponential, trigonometric and hyperbolic functions can also be expressed
with the aid of the bases introduced in Eq. (\ref{e11}).
Using the expression of the n-complex number in Eq. (\ref{e13a})
yields for the exponential of the n-complex variable $u$
\begin{eqnarray}
e^u= 
\sum_{k=1}^{n/2}e^{v_k}\left(e_k \cos \tilde v_k+\tilde e_k \sin\tilde
v_k\right).
\label{73a}
\end{eqnarray}

The trigonometric functions can be obtained from Eq. (\ref{73a}) 
with the aid of Eqs. (\ref{64}). The trigonometric functions of the
n-complex variable $u$ are
\begin{equation}
\cos u=\sum_{k=1}^{n/2}\left(e_k \cos v_k\cosh \tilde v_k
-\tilde e_k \sin v_k\sinh\tilde v_k\right),
\label{74a}
\end{equation}
\begin{equation}
\sin u= 
\sum_{k=1}^{n/2}\left(e_k \sin v_k\cosh \tilde v_k
+\tilde e_k \cos v_k\sinh\tilde v_k\right).
\label{74b}
\end{equation}

The hyperbolic functions can be obtained from Eq. (\ref{73a}) 
with the aid of Eqs. (\ref{70}). The hyperbolic functions of the
n-complex variable $u$ are
\begin{equation}
\cosh u=
\sum_{k=1}^{n/2}\left(e_k \cosh v_k\cos \tilde v_k
+\tilde e_k \sinh v_k\sin\tilde v_k\right),
\label{75a}
\end{equation}
\begin{equation}
\sinh u=
\sum_{k=1}^{n/2}\left(e_k \sinh v_k\cos \tilde v_k
+\tilde e_k \cosh v_k\sin\tilde v_k\right).
\label{75b}
\end{equation}

\section{Power series of planar n-complex numbers}

An n-complex series is an infinite sum of the form
\begin{equation}
a_0+a_1+a_2+\cdots+a_n+\cdots , 
\label{76}
\end{equation}
where the coefficients $a_n$ are n-complex numbers. The convergence of 
the series (\ref{76}) can be defined in terms of the convergence of its $n$
real components. The convergence of a n-complex series can also be studied
using n-complex variables. The main criterion for absolute convergence 
remains the comparison theorem, but this requires a number of inequalities
which will be discussed further.

The modulus $d=|u|$ of an n-complex number $u$ has been defined in Eq.
(\ref{10}). Since $|x_0|\leq |u|, |x_1|\leq |u|,..., |x_{n-1}|\leq |u|$, a
property of absolute convergence established via a comparison theorem based on
the modulus of the series (\ref{76}) will ensure the absolute convergence of
each real component of that series.

The modulus of the sum $u_1+u_2$ of the n-complex numbers $u_1, u_2$ fulfils
the inequality
\begin{equation}
||u^\prime|-|u^{\prime\prime}||\leq |u^\prime+u^{\prime\prime}|\leq
|u^\prime|+|u^{\prime\prime}| . 
\label{78}
\end{equation}
For the product, the relation is 
\begin{equation}
|u^\prime u^{\prime\prime}|\leq \sqrt{\frac{n}{2}}|u^\prime|
|u^{\prime\prime}|, 
\label{79}
\end{equation}
as can be shown from Eq. (\ref{17}). The relation (\ref{79})
replaces the relation of equality extant between 2-dimensional regular complex
numbers. 

For $u=u^\prime$ Eq. (\ref{79}) becomes
\begin{equation}
|u^2|\leq \sqrt{\frac{n}{2}} |u|^2 ,
\label{80}
\end{equation}
and in general
\begin{equation}
|u^l|\leq \left(\frac{n}{2}\right)^{(l-1)/2}|u|^l ,
\label{81}
\end{equation}
where $l$ is a natural number.
From Eqs. (\ref{79}) and (\ref{81}) it results that
\begin{equation}
|au^l|\leq \left(\frac{n}{2}\right)^{l/2} |a| |u|^l .
\label{82}
\end{equation}

A power series of the n-complex variable $u$ is a series of the form
\begin{equation}
a_0+a_1 u + a_2 u^2+\cdots +a_l u^l+\cdots .
\label{83}
\end{equation}
Since
\begin{equation}
\left|\sum_{l=0}^\infty a_l u^l\right| \leq  \sum_{l=0}^\infty
(n/2)^{l/2} |a_l| |u|^l ,
\label{84}
\end{equation}
a sufficient condition for the absolute convergence of this series is that
\begin{equation}
\lim_{l\rightarrow \infty}\frac{\sqrt{n/2}|a_{l+1}||u|}{|a_l|}<1 .
\label{85}
\end{equation}
Thus the series is absolutely convergent for 
\begin{equation}
|u|<c,
\label{86}
\end{equation}
where 
\begin{equation}
c=\lim_{l\rightarrow\infty} \frac{|a_l|}{\sqrt{n/2}|a_{l+1}|} .
\label{87}
\end{equation}

The convergence of the series (\ref{83}) can be also studied with the aid of
the formula (\ref{59a}) which is
valid for any values of $x_0,...,x_{n-1}$, as mentioned previously.
If $a_l=\sum_{p=0}^{n-1} h_p a_{lp}$, and
\begin{equation}
A_{lk}=\sum_{p=0}^{n-1} a_{lp}\cos\frac{\pi (2k-1)p}{n},
\label{88b}
\end{equation}
\begin{equation}
\tilde A_{lk}=\sum_{p=0}^{n-1} a_{lp}\sin\frac{\pi (2k-1)p}{n},
\label{88c}
\end{equation}
where $k=1,...,n/2$, the series (\ref{83}) can be written as
\begin{equation}
\sum_{l=0}^\infty \left[
\sum_{k=1}^{n/2}
(e_k A_{lk}+\tilde e_k\tilde A_{lk})(e_k v_k+\tilde e_k\tilde v_k)^l 
\right].
\label{89a}
\end{equation}
The series in Eq. (\ref{89a}) can be regarded as the sum of the $n/2$ series
obtained from each value of $k$, so that the series in Eq. (\ref{83}) is
absolutely convergent for   
\begin{equation}
\rho_k<c_k, 
\label{90}
\end{equation}
for $k=1,..., n/2$, where 
\begin{equation}
c_k=\lim_{l\rightarrow\infty} \frac
{\left[A_{lk}^2+\tilde A_{lk}^2\right]^{1/2}}
{\left[A_{l+1,k}^2+\tilde A_{l+1,k}^2\right]^{1/2}} .
\label{91}
\end{equation}
The relations (\ref{90}) show that the region of convergence of the series
(\ref{83}) is an n-dimensional cylinder.

It can be shown that $c=\sqrt{2/n}\;{\rm
min}(c_+,c_-,c_1,...,c_{n/2-1})$, where ${\rm min}$ designates the smallest of
the numbers in the argument of this function. Using the expression of $|u|$ in
Eq. (\ref{17}),  it can be seen that the spherical region of
convergence defined in Eqs. (\ref{86}), (\ref{87}) is a subset of the
cylindrical region of convergence defined in Eqs. (\ref{90}) and (\ref{91}).

\section{Analytic functions of planar n-complex variables}

The derivative  
of a function $f(u)$ of the n-complex variables $u$ is
defined as a function $f^\prime (u)$ having the property that
\begin{equation}
|f(u)-f(u_0)-f^\prime (u_0)(u-u_0)|\rightarrow 0 \:\:{\rm as} 
\:\:|u-u_0|\rightarrow 0 . 
\label{h88}
\end{equation}
If the difference $u-u_0$ is not parallel to one of the nodal hypersurfaces,
the definition in Eq. (\ref{h88}) can also 
be written as
\begin{equation}
f^\prime (u_0)=\lim_{u\rightarrow u_0}\frac{f(u)-f(u_0)}{u-u_0} .
\label{h89}
\end{equation}
The derivative of the function $f(u)=u^m $, with $m$ an integer, 
is $f^\prime (u)=mu^{m-1}$, as can be seen by developing $u^m=[u_0+(u-u_0)]^m$
as
\begin{equation}
u^m=\sum_{p=0}^{m}\frac{m!}{p!(m-p)!}u_0^{m-p}(u-u_0)^p,
\label{h90}
\end{equation}
and using the definition (\ref{h88}).

If the function $f^\prime (u)$ defined in Eq. (\ref{h88}) is independent of the
direction in space along which $u$ is approaching $u_0$, the function $f(u)$ 
is said to be analytic, analogously to the case of functions of regular complex
variables. \cite{3} 
The function $u^m$, with $m$ an integer, 
of the n-complex variable $u$ is analytic, because the
difference $u^m-u_0^m$ is always proportional to $u-u_0$, as can be seen from
Eq. (\ref{h90}). Then series of
integer powers of $u$ will also be analytic functions of the n-complex
variable $u$, and this result holds in fact for any commutative algebra. 

If an analytic function is defined by a series around a certain point, for
example $u=0$, as
\begin{equation}
f(u)=\sum_{k=0}^\infty a_k u^k ,
\label{h91a}
\end{equation}
an expansion of $f(u)$ around a different point $u_0$,
\begin{equation}
f(u)=\sum_{k=0}^\infty c_k (u-u_0)^k ,
\label{h91aa}
\end{equation}
can be obtained by
substituting in Eq. (\ref{h91a}) the expression of $u^k$ according to Eq.
(\ref{h90}). Assuming that the series are absolutely convergent so that the
order of the terms can be modified and ordering the terms in the resulting
expression according to the increasing powers of $u-u_0$ yields
\begin{equation}
f(u)=\sum_{k,l=0}^\infty \frac{(k+l)!}{k!l!}a_{k+l} u_0^l (u-u_0)^k .
\label{h91b}
\end{equation}
Since the derivative of order $k$ at $u=u_0$ of the function $f(u)$ , Eq.
(\ref{h91a}), is 
\begin{equation}
f^{(k)}(u_0)=\sum_{l=0}^\infty \frac{(k+l)!}{l!}a_{k+l} u_0^l ,
\label{h91c}
\end{equation}
the expansion of $f(u)$ around $u=u_0$, Eq. (\ref{h91b}), becomes
\begin{equation}
f(u)=\sum_{k=0}^\infty \frac{1}{k!} f^{(k)}(u_0)(u-u_0)^k ,
\label{h91d}
\end{equation}
which has the same form as the series expansion of 2-dimensional complex
functions. 
The relation (\ref{h91d}) shows that the coefficients in the series expansion,
Eq. (\ref{h91aa}), are
\begin{equation}
c_k=\frac{1}{k!}f^{(k)}(u_0) .
\label{h92}
\end{equation}

The rules for obtaining the derivatives and the integrals of the basic
functions can 
be obtained from the series of definitions and, as long as these series
expansions have the same form as the corresponding series for the
2-dimensional complex functions, the rules of derivation and integration remain
unchanged. 

If the n-complex function $f(u)$
of the n-complex variable $u$ is written in terms of 
the real functions $P_k(x_0,...,x_{n-1}), k=0,1,...,n-1$ of the real
variables $x_0,x_1,...,x_{n-1}$ as 
\begin{equation}
f(u)=\sum_{k=0}^{n-1}h_kP_k(x_0,...,x_{n-1}),
\label{h93}
\end{equation}
where $h_0=1$, then relations of equality 
exist between the partial derivatives of the functions $P_k$. 
The derivative of the function $f$ can be written as
\begin{eqnarray}
\lim_{\Delta u\rightarrow 0}\frac{1}{\Delta u} 
\sum_{k=0}^{n-1}\left(h_k\sum_{l=0}^{n-1}
\frac{\partial P_k}{\partial x_l}\Delta x_l\right),
\label{h94}
\end{eqnarray}
where
\begin{equation}
\Delta u=\sum_{k=0}^{n-1}h_l\Delta x_l.
\label{h94a}
\end{equation}
The relations between the partials derivatives of the functions $P_k$ are
obtained by setting successively in   
Eq. (\ref{h94}) $\Delta u=h_l\Delta x_l$, for $l=0,1,...,n-1$, and equating the
resulting expressions. 
The relations are 
\begin{equation}
\frac{\partial P_k}{\partial x_0} = \frac{\partial P_{k+1}}{\partial x_1} 
=\cdots=\frac{\partial P_{n-1}}{\partial x_{n-k-1}} 
= -\frac{\partial P_0}{\partial x_{n-k}}=\cdots
=-\frac{\partial P_{k-1}}{\partial x_{n-1}}, 
\label{h95}
\end{equation}
for $k=0,1,...,n-1$.
The relations (\ref{h95}) are analogous to the Riemann relations
for the real and imaginary components of a complex function. 
It can be shown from Eqs. (\ref{h95}) that the components $P_k$ fulfil the
second-order equations
\begin{eqnarray}
\lefteqn{\frac{\partial^2 P_k}{\partial x_0\partial x_l}
=\frac{\partial^2 P_k}{\partial x_1\partial x_{l-1}}
=\cdots=
\frac{\partial^2 P_k}{\partial x_{[l/2]}\partial x_{l-[l/2]}}}\nonumber\\
&&=-\frac{\partial^2 P_k}{\partial x_{l+1}\partial x_{n-1}}
=-\frac{\partial^2 P_k}{\partial x_{l+2}\partial x_{n-2}}
=\cdots
=-\frac{\partial^2 P_k}{\partial x_{l+1+[(n-l-2)/2]}
\partial x_{n-1-[(n-l-2)/2]}} ,
\label{96}
\end{eqnarray}
for $k,l=0,1,...,n-1$.

\section{Integrals of planar n-complex functions}

The singularities of n-complex functions arise from terms of the form
$1/(u-u_0)^m$, with $m>0$. Functions containing such terms are singular not
only at $u=u_0$, but also at all points of the hypersurfaces
passing through the pole $u_0$ and which are parallel to the nodal
hypersurfaces.  

The integral of an n-complex function between two points $A, B$ along a path
situated in a region free of singularities is independent of path, which means
that the integral of an analytic function along a loop situated in a region
free of singularities is zero,
\begin{equation}
\oint_\Gamma f(u) du = 0,
\label{111}
\end{equation}
where it is supposed that a surface $\Sigma$ spanning 
the closed loop $\Gamma$ is not intersected by any of
the hypersurfaces associated with the
singularities of the function $f(u)$. Using the expression, Eq. (\ref{h93}),
for $f(u)$ and the fact that 
\begin{eqnarray}
du=\sum_{k=0}^{n-1}h_k dx_k, 
\label{111a}
\end{eqnarray}
the explicit form of the integral in Eq. (\ref{111}) is
\begin{eqnarray}
\oint _\Gamma f(u) du = \oint_\Gamma
\sum_{k=0}^{n-1}h_k\sum_{l=0}^{n-1}(-1)^{[(n-k-1+l)/n]}
P_l dx_{k-l+n[(n-k-1+l)/n]}.
\label{112}
\end{eqnarray}

If the functions $P_k$ are regular on a surface $\Sigma$
spanning the loop $\Gamma$,
the integral along the loop $\Gamma$ can be transformed in an integral over the
surface $\Sigma$ of terms of the form
$\partial P_l/\partial x_{k-m+n[(n-k+m-1)/n]} 
- (-1)^s \partial P_m/\partial x_{k-l+n[(n-k+l-1)/n]}$, where 
$s=[(n-k+m-1)/n]-[(n-k+l-1)/n]$.
These terms are equal to zero by Eqs. (\ref{h95}), and this
proves Eq. (\ref{111}). 

The integral of the function $(u-u_0)^m$ on a closed loop $\Gamma$ is equal to
zero for $m$ a positive or negative integer not equal to -1,
\begin{equation}
\oint_\Gamma (u-u_0)^m du = 0, \:\: m \:\:{\rm integer},\: m\not=-1 .
\label{112b}
\end{equation}
This is due to the fact that $\int (u-u_0)^m du=(u-u_0)^{m+1}/(m+1), $ and to
the fact that the function $(u-u_0)^{m+1}$ is singlevalued for $n$ an integer.

The integral $\oint_\Gamma du/(u-u_0)$ can be calculated using the exponential
form, Eq. (\ref{50a}), for the difference $u-u_0$, 
\begin{eqnarray}
u-u_0=\rho\exp\left\{\sum_{p=1}^{n-1}h_p\left[
-\frac{2}{n}\sum_{k=2}^{n/2}
\cos\left(\frac{2\pi kp}{n}\right)\ln\tan\psi_{k-1}
\right]
+\sum_{k=1}^{n/2}\tilde e_k\phi_k\right\}.
\label{113a}
\end{eqnarray}
Thus the quantity $du/(u-u_0)$ is
\begin{eqnarray}
\frac{du}{u-u_0}=
\frac{d\rho}{\rho}+\sum_{p=1}^{n-1}h_p\left[
-\frac{2}{n}\sum_{k=2}^{n/2}
\cos\left(\frac{2\pi kp}{n}\right)d\ln\tan\psi_{k-1}\right]
+\sum_{k=1}^{n/2}\tilde e_kd\phi_k.
\label{114a}
\end{eqnarray}
Since $\rho$ and $\ln(\tan\psi_{k-1})$ are singlevalued variables, it follows
that 
$\oint_\Gamma d\rho/\rho =0$, and 
$\oint_\Gamma d(\ln\tan\psi_{k-1})=0$.
On the other hand, since $\phi_k$ are cyclic variables, they may give
contributions to the integral around the closed loop $\Gamma$.

The expression of $\oint_\Gamma du/(u-u_0)$ can be written 
with the aid of a functional which will be called int($M,C$), defined for a
point $M$ and a closed curve $C$ in a two-dimensional plane, such that 
\begin{equation}
{\rm int}(M,C)=\left\{
\begin{array}{l}
1 \;\:{\rm if} \;\:M \;\:{\rm is \;\:an \;\:interior \;\:point \;\:of} \;\:C
,\\  
0 \;\:{\rm if} \;\:M \;\:{\rm is \;\:exterior \;\:to}\:\; C .\\
\end{array}\right.
\label{118}
\end{equation}
With this notation the result of the integration on a closed path $\Gamma$
can be written as 
\begin{equation}
\oint_\Gamma\frac{du}{u-u_0}=
\sum_{k=1}^{n/2}2\pi\tilde e_k 
\;{\rm int}(u_{0\xi_k\eta_k},\Gamma_{\xi_k\eta_k}) ,
\label{119}
\end{equation}
where $u_{0\xi_k\eta_k}$ and $\Gamma_{\xi_k\eta_k}$ are respectively the
projections of the point $u_0$ and of 
the loop $\Gamma$ on the plane defined by the axes $\xi_k$ and $\eta_k$,
as shown in Fig. 3.

If $f(u)$ is an analytic n-complex function which can be expanded in a
series as written in Eq. (\ref{h91aa}), and the expansion holds on the curve
$\Gamma$ and on a surface spanning $\Gamma$, then from Eqs. (\ref{112b}) and
(\ref{119}) it follows that
\begin{equation}
\oint_\Gamma \frac{f(u)du}{u-u_0}=
2\pi f(u_0)\sum_{k=1}^{n/2}\tilde e_k 
\;{\rm int}(u_{0\xi_k\eta_k},\Gamma_{\xi_k\eta_k}) .
\label{120}
\end{equation}

Substituting in the right-hand side of 
Eq. (\ref{120}) the expression of $f(u)$ in terms of the real 
components $P_k$, Eq. (\ref{h93}), yields
\begin{eqnarray}
\lefteqn{\oint_\Gamma \frac{f(u)du}{u-u_0}
=\frac{2}{n}\sum_{k=1}^{n/2}\sum_{l=0}^{n-1}h_l
\nonumber}\\
&&\sum_{p=1}^{n-1}
(-1)^{[(l-p)/n]}\sin\left[\frac{\pi (2k-1)p}{n}\right] 
P_{n-p+l-n[(n-p+l)/n]}(u_0)
\:{\rm int}(u_{0\xi_k\eta_k},\Gamma_{\xi_k\eta_k}) .
\label{121}
\end{eqnarray}
It the integral in Eq. (\ref{121}) is written as 
\begin{equation}
\oint_\Gamma \frac{f(u)du}{u-u_0}=\sum_{l=0}^{n-1}h_l I_l,
\label{122a}
\end{equation}
it can be checked that
\begin{equation}
\sum_{l=0}^{n-1} I_l=0.
\label{122b}
\end{equation}

If $f(u)$ can be expanded as written in Eq. (\ref{h91aa}) on 
$\Gamma$ and on a surface spanning $\Gamma$, then from Eqs. (\ref{112b}) and
(\ref{119}) it also results that
\begin{equation}
\oint_\Gamma \frac{f(u)du}{(u-u_0)^{n+1}}=
\frac{2\pi}{n!}f^{(n)}(u_0)\sum_{k=1}^{[(n-1)/2]}\tilde e_k 
\;{\rm int}(u_{0\xi_k\eta_k},\Gamma_{\xi_k\eta_k}) ,
\label{122}
\end{equation}
where the fact has been used  that the derivative $f^{(n)}(u_0)$ is related to
the expansion coefficient in Eq. (\ref{h91aa}) according to Eq. (\ref{h92}).

If a function $f(u)$ is expanded in positive and negative powers of $u-u_l$,
where $u_l$ are n-complex constants, $l$ being an index, the integral of $f$
on a closed loop $\Gamma$ is determined by the terms in the expansion of $f$
which are of the form $r_l/(u-u_l)$,
\begin{equation}
f(u)=\cdots+\sum_l\frac{r_l}{u-u_l}+\cdots.
\label{123}
\end{equation}
Then the integral of $f$ on a closed loop $\Gamma$ is
\begin{equation}
\oint_\Gamma f(u) du = 
2\pi \sum_l\sum_{k=1}^{n/2}\tilde e_k 
\;{\rm int}(u_{l\xi_k\eta_k},\Gamma_{\xi_k\eta_k})r_l .
\label{124}
\end{equation}

\section{Factorization of planar n-complex polynomials}

A polynomial of degree $m$ of the n-complex variable $u$ has the form
\begin{equation}
P_m(u)=u^m+a_1 u^{m-1}+\cdots+a_{m-1} u +a_m ,
\label{125}
\end{equation}
where $a_l$, for $l=1,...,m$, are in general n-complex constants.
If $a_l=\sum_{p=0}^{n-1} h_p a_{lp}$, and with the
notations of Eqs. (\ref{88b}), (\ref{88c}) applied for $l= 1, \cdots, m$, the
polynomial $P_m(u)$ can be written as
\begin{eqnarray}
P_m= \sum_{k=1}^{n/2}\left\{(e_k v_k+\tilde e_k\tilde v_k)^m+
\sum_{l=1}^m(e_k A_{lk}+\tilde e_k\tilde A_{lk})
(e_k v_k+\tilde e_k\tilde v_k)^{m-l} 
\right\},
\label{126a}
\end{eqnarray}
where the constants $A_{lk}, \tilde A_{lk}$ are real numbers.

The polynomials of degree $m$ in $e_k v_k+\tilde e_k\tilde v_k$ 
in Eq. (\ref{126a}) 
can always be written as a product of linear factors of the form
$e_k (v_k-v_{kp})+\tilde e_k(\tilde v_k-\tilde v_{kp})$, where the
constants $v_{kp}, \tilde v_{kp}$ are real,
\begin{eqnarray}
\lefteqn{(e_k v_k+\tilde e_k\tilde v_k)^m+
\sum_{l=1}^m(e_k A_{lk}+\tilde e_k\tilde A_{lk})
(e_k v_k+\tilde e_k\tilde v_k)^{m-l} \nonumber}\\
&&=\prod_{p=1}^{m}\left\{e_k (v_k-v_{kp})
+\tilde e_k(\tilde v_k-\tilde v_{kp})\right\}.
\label{126b}
\end{eqnarray}
Then the polynomial $P_m$ can be written as
\begin{equation}
P_m=\sum_{k=1}^{n/2}\prod_{p=1}^m
\left\{e_k (v_k-v_{kp})+\tilde e_k(\tilde v_k-\tilde v_{kp})\right\}.
\label{127a}
\end{equation}
Due to the relations  (\ref{e12a}),
the polynomial $P_m(u)$ can be written as a product of factors of
the form  
\begin{eqnarray}
P_m(u)=\prod_{p=1}^m \left\{\sum_{k=1}^{n/2}
\left\{e_k (v_k-v_{kp})
+\tilde e_k(\tilde v_k-\tilde v_{kp})\right\}\right\}.
\label{128}
\end{eqnarray}
This relation can be written with the aid of Eq. (\ref{e13a}) as
\begin{eqnarray}
P_m(u)=\prod_{p=1}^m (u-u_p) ,
\label{129a}
\end{eqnarray}
where
\begin{eqnarray}
u_p=\sum_{k=1}^{n/2}\left(e_k v_{kp}
+\tilde e_k\tilde v_{kp}\right), 
\label{129b}
\end{eqnarray}
for $p=1,...,m$.
For a given $k$, the roots  
$e_k v_{k1}+\tilde e_k\tilde v_{k1}, ...,  e_k v_{km}+\tilde e_k\tilde v_{km}$
defined in Eq. (\ref{126b}) may be ordered arbitrarily. This means that Eq.
(\ref{129b}) gives sets of $m$ roots 
$u_1,...,u_m$ of the polynomial $P_m(u)$, 
corresponding to the various ways in which the roots $e_k v_{kp}+\tilde
e_k\tilde v_{kp}$ are ordered according to $p$ for each value of $k$. 
Thus, while the n-complex components in Eq. (\ref{126b}) taken separately
have  
unique factorizations, the polynomial $P_m(u)$ can be written in many different
ways as a product of linear factors. 

If $P(u)=u^2+1$, the degree is $m=2$, the coefficients of the polynomial are
$a_1=0, a_2=1$, the n-complex components of $a_2$ are $a_{20}=1, a_{21}=0,
... ,a_{2,n-1}=0$, the components $A_{2k}, \tilde A_{2k}$ calculated according
to Eqs. (\ref{88b}), (\ref{88c}) are $A_{2k}=1, \tilde A_{2k}=0, k=1,...,n/2$.
The left-hand side of Eq. (\ref{126b}) has the form 
$(e_k v_k+\tilde e_k\tilde v_k)^2+e_k$, and since $e_k=-\tilde e_k^2$, the
right-hand side of Eq. (\ref{126b}) is 
$\left\{e_k v_k+\tilde e_k(\tilde v_k+1)\right\}
\left\{e_k v_k+\tilde e_k(\tilde v_k-1)\right\}$, so that
$v_{kp}=0, \tilde v_{kp}=\pm 1, k=1,...,n/2, p=1,2$.
Then Eq. (\ref{127a}) has the form $u^2+1=\sum_{k=1}^{n/2}
\left\{e_k v_k+\tilde e_k(\tilde v_k+1)\right\}
\left\{e_k v_k+\tilde e_k(\tilde v_k-1)\right\}$.
The factorization in Eq. (\ref{129a}) is $u^2+1=(u-u_1)(u-u_2)$, where 
$u_1=\pm \tilde e_1\pm\tilde e_2\pm\cdots\pm \tilde e_{n/2}, u_2=-u_1$, so that
there are $2^{n/2-1}$ independent sets of roots $u_1,u_2$
of $u^2+1$. It can be checked
that $(\pm \tilde e_1\pm\tilde e_2\pm\cdots\pm \tilde e_{n/2})^2=
-e_1-e_2-\cdots-e_{n/2}=-1$.

\section{Representation of planar n-complex numbers by irreducible matrices}

If the unitary matrix written in Eq. (\ref{11}) is called $T$,
it can be shown that the matrix $T U T^{-1}$ has the form 
\begin{equation}
T U T^{-1}=\left(
\begin{array}{cccc}
V_1      &     0   & \cdots  &   0   \\
0        &     V_2 & \cdots  &   0   \\
\vdots   &  \vdots & \cdots  &\vdots \\
0        &     0   & \cdots  &   V_{n/2}\\
\end{array}
\right),
\label{129}
\end{equation}
where $U$ is the matrix in Eq. (\ref{24a}) used to represent the n-complex
number $u$. In Eq. (\ref{129}), the matrices $V_k$ are
the matrices
\begin{equation}
V_k=\left(
\begin{array}{cc}
v_k           &     \tilde v_k   \\
-\tilde v_k   &     v_k          \\
\end{array}\right),
\label{130}
\end{equation}
for $ k=1,...,n/2$, where $v_k, \tilde v_k$ are the variables introduced in Eqs. (\ref{9a}) and
(\ref{9b}), and the symbols 0 denote
the matrix
\begin{equation}
\left(
\begin{array}{cc}
0   &  0   \\
0   &  0   \\
\end{array}\right).
\label{131}
\end{equation}
The relations between the variables $v_k, \tilde v_k$ for the multiplication of
n-complex numbers have been written in Eq. (\ref{22}). The matrix
$T U T^{-1}$ provides an irreducible representation
\cite{4} of the n-complex number $u$ in terms of matrices with real
coefficients. For $n=2$, Eqs. (\ref{9a}) and (\ref{9b}) give $v_1=x_0, 
\tilde v_1=x_1$, and Eq. (\ref{e11}) gives $e_1=1, \tilde e_1=h_1$, 
where according to Eq. (\ref{1}) $h_1^2=-1$, so that the matrix $V_1$, Eq.
(\ref{130}), is 
\begin{equation}
v_1=\left(
\begin{array}{cc}
x_0    &  x_1   \\
-x_1   &  x_0   \\
\end{array}\right),
\label{132}
\end{equation}
which shows that, for $n=2$, the hypercomplex numbers $x_0+h_1 x_1$ are
identical to the usual 2-dimensional complex numbers $x+iy$. 

\section{Conclusions}

The operations of addition and multiplication of the n-complex numbers
introduced in this 
work have a geometric interpretation based on the amplitude $\rho$,
the modulus $d$, the planar angles $\psi_k$ and the azimuthal angles $\phi_k$. 
The n-complex numbers can be written in exponential and
trigonometric forms with the aid of these variables.
The n-complex functions defined by series of powers are analytic, and 
the partial derivatives of the components of the n-complex functions are
closely related. The integrals of n-complex functions are independent of path
in regions where the functions are regular. The fact that the exponential form
of the n-complex numbers depends on the cyclic variables $\phi_k$
leads to the 
concept of pole and residue for integrals on closed paths. The polynomials of
n-complex variables can always be written as products of linear factors,
although the factorization is not unique.

\newpage

FIGURE CAPTIONS\\

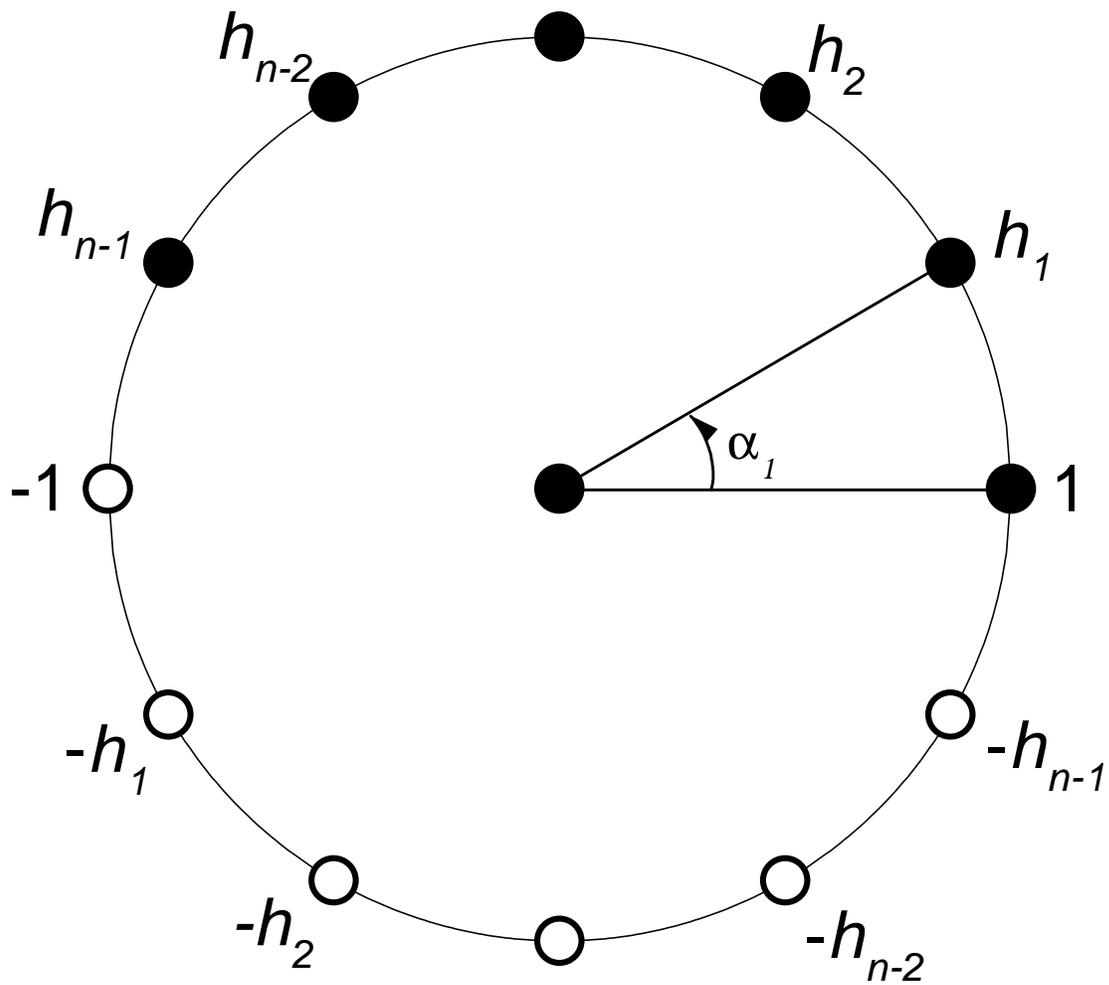
Fig. 1. Representation of the hypercomplex bases $1, h_1,...,h_{n-1}$
by points on a circle at the angles $\alpha_k=\pi k/n$.
The product $h_j h_k$ will be represented by the point of the circle at the
angle $\pi (j+k)/2n$, $j,k=0,1,...,n-1$. If $\pi\leq\pi (j+k)/2n\leq 2\pi$, the
point is opposite to the basis $h_l$ of angle $\alpha_l=\pi (j+k)/n-\pi$. \\

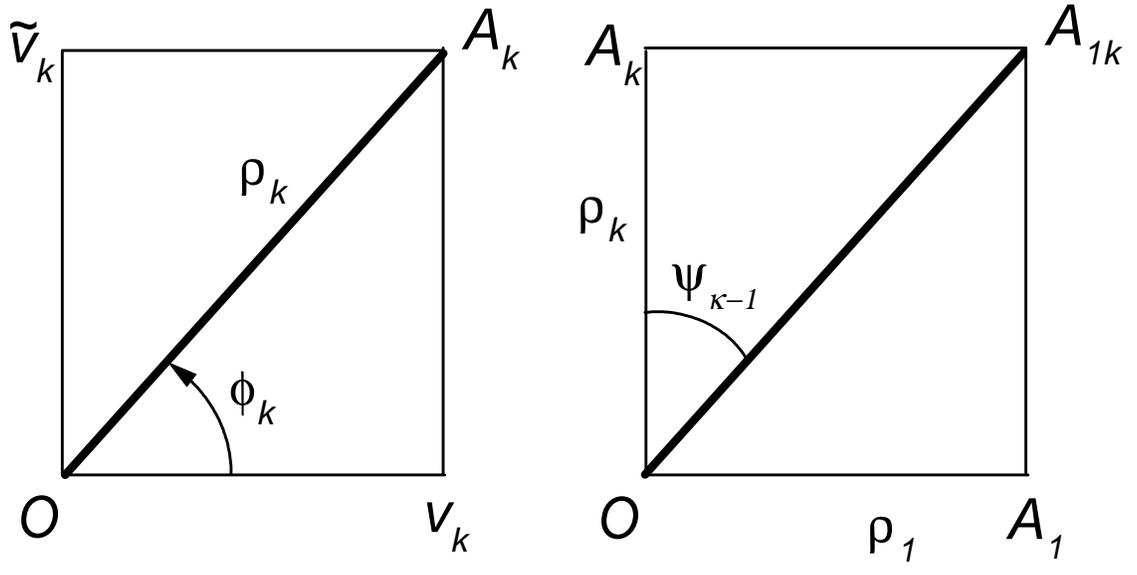
Fig. 2. Radial distance $\rho_k$  and
azimuthal angle $\phi_k$ in the plane of the axes $v_k,\tilde v_k$, and planar
angle $\psi_{k-1}$ between the line $OA_{1k}$ and the 2-dimensional plane
defined by the axes $v_k,\tilde v_k$. $A_k$ is the projection of the point $A$
on the plane of the axes $v_k,\tilde v_k$, and $A_{1k}$ is the projection of
the point $A$ on the 4-dimensional space defined 
by the axes $v_1, \tilde v_1, v_k,\tilde v_k$. \\

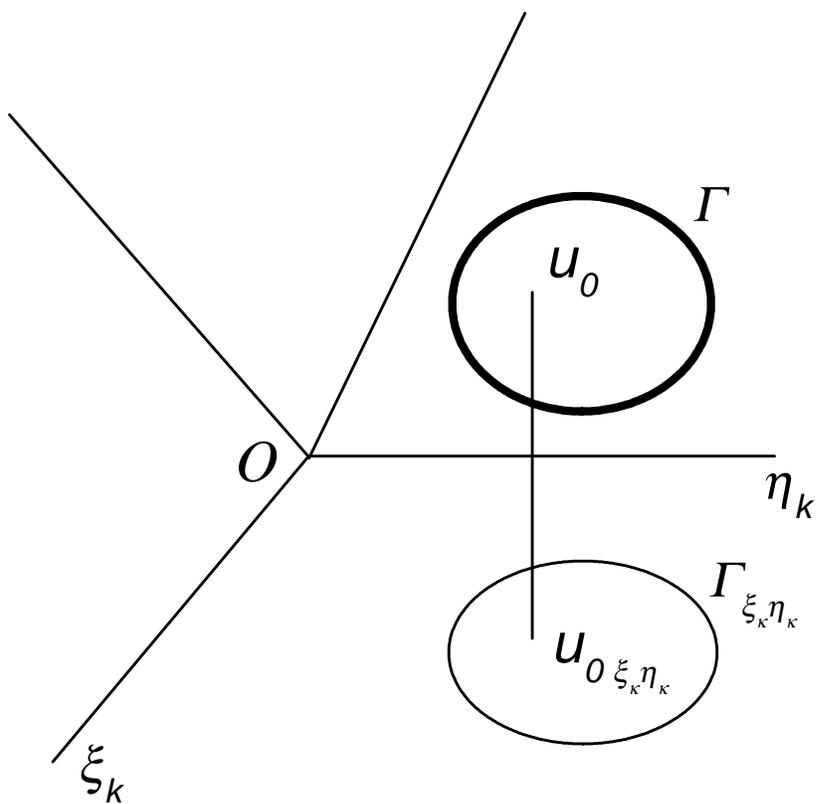
Fig. 3. Integration path $\Gamma$ and pole $u_0$, and their projections
$\Gamma_{\xi_k\eta_k}$ and $u_{0\xi_k\eta_k}$ on the plane $\xi_k \eta_k$.\\ 

\end{document}